\documentclass[reqno]{amsart}

\usepackage{setspace,enumitem,amsmath,amssymb,amsfonts,color,graphicx,bigints,hyperref}
\usepackage{tikz,array} %This package offers the ability to draw pictures
\usepackage{tikz-cd}

\usepackage[english]{babel}
\usepackage{amsthm}

\usepackage[backend=biber,style=alphabetic, citestyle=alphabetic, sorting = nty]{biblatex}
\theoremstyle{definition}
\newtheorem{theorem}{Theorem}[section]
\newtheorem{lemma}[theorem]{Lemma}
\newtheorem{corollary}[theorem]{Corollary}

\theoremstyle{definition}
\newtheorem{definition}[theorem]{Definition}
\newtheorem{remark}[theorem]{Remark}

\title[GK Dimension and Generator Bounds]{On GK Dimension and Generator Bounds for a Class of Graded Algebras}
\author{Abdourrahmane Kabbaj}
\date{\today}
\address{Department of Mathematics, University of California, Irvine, CA 92697, USA}
\email{akabbaj@uci.edu}

\begin{document}

\maketitle

\begin{abstract} In this paper, we introduce the concept of \textit{monotonic algebras}, a broad class of algebras that includes all Artin-Schelter regular algebras of dimension at most four, as well as algebras with \textit{pure} resolutions, such as Koszul and piecewise Koszul algebras. We show that the Gelfand-Kirillov (GK) dimension of these algebras is bounded above by their global dimension and establish a similar result for the minimal number of generators. Furthermore, we prove a parity theorem for Artin-Schelter regular algebras, demonstrating that the difference between their global dimension and GK dimension is always an even integer.

\end{abstract}

\section{Introduction}

%To do: 

%Convert the Examples into a theorem, potentially expand more.

%Try to generalize results to a more general scope, maybe graded algebras of finite Global and GK dimension.

%Note: \nu(pq) \leq \nu(p) + \nu(q) is FALSE and so defining a class of algebras satisfying  $\nu(p)\leq d$ wont necessarily be closed under tensor products.
%Counter Example: (1-t)(1+t^2+t^3) = 1-t +t^2 - t^4

%Does A having finite global and GK dimension imply that h_A = q(t)/p(t) with lemma 2.7 from manny's paper applicable?

%Does Manny know of an example of an algebra of finite global dimension but one of the P_i is infinite?

Artin-Schelter regular algebras are graded algebras, often regarded as the noncommutative analogues of polynomial rings. They are algebras of finite global and GK-dimension, and their introduction marked a pivotal moment in the study of noncommutative projective geometry, where noncommutative projective spaces are constructed as the noncommutative projective scheme Proj$A$, with $A$ being an Artin-Schelter regular algebra. Introduced by Artin and Schelter in \cite{Artin1987}, they classified all regular algebras $A$ of global dimension three, generated in degree one, and showed that the minimal projective resolution of the trivial $A$-module $\mathbb{K}$ takes one of two forms:

\[   A( -3) \longrightarrow  A( -2)^3   \longrightarrow  A( -1)^3  \longrightarrow  A  \longrightarrow  \mathbb{K}    \tag{1.0.1}\label{1.0.1}  \]

\[   A( -4) \longrightarrow  A( -3)^2   \longrightarrow  A( -1)^2  \longrightarrow  A  \longrightarrow  \mathbb{K}   \tag{1.0.2}\label{1.0.2} \]\\

The algebras of global dimension three have since been fully classified by Artin, Schelter, Tate, and Van den Bergh in \cite{Artin1987}, \cite{Artin1990}, and \cite{Artin1991}, and are considered well understood by many.\\

In \cite{lu2004regular}, the authors initiated the next step in classifying Artin-Schelter regular algebras of global dimension four generated in degree one. Under the additional natural assumption that $A$ is a domain, they demonstrated that the minimal projective resolution of the trivial $A$-module $\mathbb{K}$ falls into one of three types:

 \[  A(-7) \longrightarrow A( -6)^2 \longrightarrow  A( -4)\oplus A(-3)   \longrightarrow  A( -1)^2  \longrightarrow  A  \longrightarrow  \mathbb{K}   \tag{1.1.1}\label{1.1.1} \]

 \[   A(-5) \longrightarrow A( -4)^3 \longrightarrow  A( -3)^2\oplus A(-2)^2   \longrightarrow  A( -1)^3  \longrightarrow  A  \longrightarrow  \mathbb{K}   \tag{1.1.2}\label{1.1.2} \]
 
  \[   A(-4) \longrightarrow A( -3)^4 \longrightarrow  A(-2)^6   \longrightarrow  A( -1)^4  \longrightarrow  A  \longrightarrow  \mathbb{K}  \tag{1.1.3}\label{1.1.3} \]
 They proceeded to classify those of type (\ref{1.1.1}) under some additional conditions. Algebras of type (\ref{1.1.2}) were subsequently classified under specific constraints in later work \cite{Rogalski2011}. Although a comprehensive classification of Artin-Schelter regular algebras of dimension four remains challenging, significant progress has been made in studying various families and classes (see, for example, \cite{Smith1992}, \cite{Cassidy1999}, \cite{Zhang2007DoubleER}).\\

Along with these developments, the investigation of Artin-Schelter regular algebras and their properties continues to attract considerable scholarly attention, though several fundamental questions remain unanswered. For instance, let $A$ be an Artin-Schelter regular algebra of global dimension $d$:

\begin{itemize}

\item[(1)] Is $A$ a noetherian domain?\\

\item [(2)]Does $A$ have the Hilbert series of a weighted polynomial ring. That is, $h_A(t) = \frac{1}{\prod_{i=1}^m (1-t^i)^{n_i}}$ where $n_i$'s are positive integers?\\

\item[(3)] Is GK$\dim(A) = d$?\\

\item[(4)] Is $A$ minimally generated as an algebra by  at most $d$ elements?\\

\end{itemize}

The properties  above are widely conjectured to be true by many researchers, as they have been shown to hold in all known examples. For a more detailed exploration of Artin-Schelter regular algebras, the reader is referred to the recent survey \cite{rogalski2023artinschelter}.\\

In this short paper, we advance the study of questions (3) and (4) by establishing upper bounds on the GK dimension and the number of generators of a new class of graded algebras that encompasses many classes of Artin-Schelter regular algebras that have been extensively studied. The paper is organized as follows. In Section 2, we recall the definition of Artin-Schelter regular algebras in addition to some well known results regarding their projective minimal resolutions and Hilbert series. We end the section by proving a parity result between the global and GK dimension of Artin-Schelter regular algebras. More precisely, we prove the following:

\theorem\label{theorem1.1}(Theorem \ref{theorem2.4}) Let \( A \) be an Artin-Schelter regular algebra of global dimension \( d \). Then $\mathrm{GKdim}(A) \equiv d \pmod{2}$.\\

In Section 3, we introduce the class of monotonic algebras, defined by the property that their trivial \( A \)-module, \( \mathbb{K}_A \), admits a finite free resolution with monotonically increasing Betti numbers. We establish upper bounds on their Gelfand-Kirillov dimension and under mild assumptions, derive similar bounds on the minimal number of generators. Furthermore, we show that this class includes a broad range of algebras of significant interest in the literature. Specifically, we prove the following:

\theorem\label{theorem1.2}(Theorem \ref{theorem3.4}) The class of monotonic algebras includes:\begin{enumerate}

\item[(a)] Algebras with \textit{pure} minimal resolutions such as Koszul, N-Koszul or Piecewise Koszul algebras.\\

\item[(b)] Artin-Schelter regular algebras of dimension at most four.\\

\item[(c)] Quadratic Artin-Schelter regular algebras of dimension at most six.\\

\end{enumerate}

\theorem\label{theorem1.3}(Theorem \ref{theorem3.6}) Let $A$ be a monotonic algebra of global dimension $d$. If $A$ has finite GK dimension, then GK$\dim(A) \leq  d$.

\theorem\label{theorem1.4}(Theorem \ref{theorem3.9} \&Theorem \ref{theorem3.10}  ) Let $A$ be a monotonic algebra of global dimension $d$. If either of the two conditions hold, 

\begin{itemize}

\item GK$\dim(A) = d$\\

\item $h_A(t) = \frac{1}{\prod_{i=1}^m (1-t^i)^{n_i}}$

\end{itemize}
then $A$ is generated by at most $d$ elements.

\section{Artin-Schelter Regular Algebras }

Fix a field $\mathbb{K}$.  All vector spaces will be $\mathbb{K}$-vector spaces. An algebra $A$ is called  \textit{$\mathbb N$-graded} if $A = \bigoplus_{n=0}^\infty A_n$ such that $A_iA_j \subseteq A_{i+j}$ for all $i,j$. We say  $A$ is  \textit{locally finite} if $\dim_\mathbb K A_n <\infty$ for all $n$ and \textit{connected} if $A_0 = \mathbb{K}$. Lastly, we say that  $A$ is  \textit{generated in degree one} if it has a generating set in $A_1$. All algebras $A$ in this paper will be assumed to be  $\mathbb N$-graded, locally finite, connected,  and generated  in degree one. \\

\definition\label{def2.1} A connected graded algebra $A = \bigoplus_{n=0}^\infty A_n$ is called an \textit{Artin-Schelter regular}  algebra of dimension $d$, AS regular for short, if the following conditions hold:

\begin{itemize}

\item[(1)] $A$ has finite global dimension $d<\infty$;

\item[(2)] $A$ is Gorenstein; that is
         \begin{equation*}
   \text{ \underline{Ext}}_A^i(\mathbb{K},A)=
    \begin{cases}
        0 & \text{if } i\neq d\\
        \mathbb K(l) & \text{if } i=d
    \end{cases}
    \end{equation*}
    where $\mathbb K(l)$ denotes the trivial $A$-module $\mathbb K$ in degree $-l$. The number $l$ is  called the \textit{AS-index} of $A$. 
  
   \item[(3)] It has finite Gelfand-Kirillov (GK) dimension, i.e., there exist positive constants $c$ and $d$ such that $\dim A_n \leq cn^d$ for all $n$.
   
     \end{itemize}     
      More generally, $A$ is called \textit{Artin-Schelter  Gorenstein} if it satisfies condition (2) from above and (1) is replaced by the weaker condition that $A$ has finite injective dimension $d<\infty$ as both a left and right $A$-module.\\

The following two lemmas are well known in the literature regarding the relationship between the Hilbert series, Gelfand-Kirillov (GK) dimension, and minimal projective resolutions of Artin-Schelter regular algebras, with proofs available in \cite{lu2004regular} and \cite{Artin1991}.

\lemma\label{lemma2.2} Let  $A$ be an AS Gorenstein algebra of global dimension $d$, then the following holds: 

\begin{enumerate}

\item[(i)]$A$ is finitely generated

\item[(ii)] The  trivial $A$-module $\mathbb K_A$ has a minimal free resolution of the form \[ 0\longrightarrow P_d \longrightarrow \cdots \longrightarrow P_1 \longrightarrow P_0 \longrightarrow \mathbb K_A \longrightarrow 0  \tag{2.1.1}               \]
where $P_i= \bigoplus_{s=1}^{n_i}A(-\alpha_{i,s})$ for some integers $n_i>0$ and $\alpha_{i,s}>0$.

\item[(iii)] The above resolution is symmetric in the following sense: $P_0 = A$, $P_d = A(-l)$, $n_i = n_{d-i}$, and $\alpha_{i,s} + \alpha_{d-i,n_i -s+1} =l$ for all $i,s$.\\
\end{enumerate}

\lemma\label{lemma2.3} Let  $A$ be an AS regular algebra of dimension $d$, then 

\begin{enumerate}

\item[(i)] The Hilbert series of $A$ has the form $h_A(t) = \frac{1}{p(t)}$ where \[p(t) = \sum_{i=0}^d(-1)^i \sum_{s=0}^{n_i}t^{\alpha_{i,s}}\]
with $\alpha_{s,i}$'s and $n_i$ being the same ones arising from the  minimal  projective resolution of $\mathbb K_A$ in (2.1.1).

\item[(ii)] $p(t)$ satisfies the following equality. \[      t^lp(t^{-1}) = (-1)^dp(t)     \]

\item[(iii)] The GK dimension of  $A$ is equal to the maximum order of the pole of $h_A(t)$ at $t=1$. I.e. if \[ h_A(t) =\frac{1}{(1-t)^mg(t)}               \]
and $g(1) \neq 0$, then GK$\dim A = m$.

\item[(iv)] All the roots of $p(t)$ are roots of unity. \\\\

\end{enumerate}

Let \( A \) be an Artin-Schelter regular algebra of global dimension \( d \) and Gelfand-Kirillov dimension \( m \). Furthermore, let \( \Phi_n(t) \) denote the \( n \)-th cyclotomic polynomial. By Lemma \ref{lemma2.3}, the Hilbert series of \( A \) can be expressed as  $h_A(t) = \frac{1}{p(t)}$, where  
\[
p(t) = (1-t)^m \prod_{j=1}^r \Phi_{n_j}(t)
\]  
is a product of cyclotomic polynomials with \( n_j > 1 \) for all \( j \), and \(\deg p(t) = l\), where \( l \) is the AS-index of \( A \). Because \( p(t) \) fully determines the GK dimension of \( A \), analyzing its structure can provide valuable insights into the behavior of \( \mathrm{GKdim}(A) \).\\

For instance, when working with AS-regular algebras, it is common to impose additional assumptions on the Gelfand-Kirillov dimension, such as requiring that \( \mathrm{GKdim}(A) \geq n \) for some positive integer \( n \). Examples of this approach appear in (\cite{floystad2010artinschelter}, Theorem 5.6), (\cite{Li2019}, Theorem 4.1), and (\cite{lu2004regular}, Proposition 1.4). It is notable that in these cases, the assumption  \( \mathrm{GKdim}(A) \geq n \)  frequently implies the stronger result \( \mathrm{GKdim}(A) = n+1 \).\\

In our first theorem, we leverage the fact that \( p(t) \) is a product of predominantly self-reciprocal polynomials to establish a parity relationship between the global dimension and the GK dimension of \( A \). This result provides a structural explanation for the observed phenomenon.\\

\theorem\label{theorem2.4} Let \( A \) be an Artin-Schelter regular algebra of global dimension \( d \). Then $\mathrm{GKdim}(A) \equiv d \pmod{2}$.\\

\begin{proof}Suppose that GK$\dim A =m$. The Hilbert series of $A$ has the form $h_A(t) = \frac{1}{p(t)}$ where  

\[ p(t) = (1-t)^m \prod_{j=1}^r\Phi_{n_j}(t)  \]
such that $n_j>1$ for all $j$. By Lemma \ref{lemma2.3}(ii), we also know that \[t^lp(t^{-1}) = (-1)^dp(t)\tag{2.5.1}\label{2.5.1}\]
Expanding the left hand side of (\ref{2.5.1}), we get 
\[  t^lp(t^{-1}) =  t^l(1-t^{-1})^m\prod_{j=1}^r\Phi_{n_j}(t^{-1})       \]
\[  = t^l(1-t^{-1})^m \prod_{j=1}^rt^{-\varphi(\alpha_j)}\cdot\Phi_{n_j}(t)    \]
\[  = (-1)^m(1-t)^m  \prod_{j=1}^r\Phi_{n_j}(t)              \]
where the second equality follows from the well known fact that for $n>1$, cyclotomic polynomials $\Phi_n(x)$ are self reciprocal. \\
Now expanding the right hand side of (\ref{2.5.1}), we get 
 \[    (-1)^dp(t) = (-1)^d(1-t)^m \prod_{j=1}^r\Phi_{n_j}(t) \]
Since both sides are equal, we conclude that $(-1)^d = (-1)^m$ and thus  either both $d$ and $m$ are even, or both are odd.

\end{proof}

This theorem can offer a significant simplification in reducing manual computations. For instance, in (\cite{Li2019}, Theorem 4.1), the author proves that if $A$ is a Koszul Artin-Schelter regular algebra of dimension 5, with $\mathrm{GKdim}(A) \geq 4$, then $h_A(t) = \frac{1}{(1-t)^5}$. The proof involves constructing a generic minimal resolution, computing multiple derivatives, and performing extensive algebraic manipulations. As we show below, such a statement and its generalization is more or less immediate by Theorem 2.4.\\

\corollary\label{corollary2.5}  
Let \( A \) be a Koszul Artin-Schelter regular algebra of global dimension \( d \). Suppose further that \( \mathrm{GKdim}(A) \geq d-1 \). Then,  
\[
h_A(t) = \frac{1}{(1-t)^d}.
\]

\begin{proof}  
Since \( A \) is Koszul AS-regular, its AS-index satisfies \( l = d \), and its Hilbert series takes the form \( h_A(t) = \frac{1}{p(t)} \), where \( \deg p(t) = d \).  By Lemma \ref{lemma2.3}(iii), this immediately implies that \( \mathrm{GKdim}(A) \leq d \). Given the assumption that \( \mathrm{GKdim}(A) \geq d -1 \), applying Theorem \ref{theorem2.4} yields \( \mathrm{GKdim}(A) = d \) and the result follows.\\

\end{proof}

\section{Monotonic Algebras of Finite Global Dimension}

Consider a connected graded algebra $A$ with finite global dimension $d$. Since every graded left-bounded projective $A$-module is free, it can be expressed as a sum of shifts of $A$. Thus, any left-bounded $A$-module $M$ admits a minimal free resolution of the form 

\[ 0\longrightarrow P_d \longrightarrow \cdots \longrightarrow P_1 \longrightarrow P_0 \longrightarrow \mathbb K_A \longrightarrow 0  \tag{3.0.1}\label{3.0.1}               \]
where $P_i= \bigoplus_{s=1}^{n_i}A(-\alpha_{i,s})$, and $n_i$ are possibly infinite. If $n_i$ is finite for all $i$, $M$ is  said to posses a \textit{finite free resolution}. Following the terminology in \cite{Stephenson1997}, if $M$ has a finite free resolution, its \textit{characteristic polynomial} is defined as 
\[c_M (t) = \sum_{i=0}^d(-1)^i \sum_{s=0}^{n_i}t^{\alpha_{i,s}}.\tag{3.0.2}\label{3.0.2}\]\\

Suppose that \( A \) is a graded algebra of global dimension \( d \) for which the trivial \( A \)-module \( \mathbb{K}_A \) admits a finite free resolution as in (\ref{3.0.1}). Let \( p(t) \) denote the characteristic polynomial of \( \mathbb{K}_A \). By Lemma 2.3 of \cite{Stephenson1997}, the Hilbert series of \( A \) is given by  $h_A(t) = \frac{1}{p(t)}$.\\

The graded Betti numbers ($\alpha_{i,s}$),  appearing in $p(t)$ are crucial in defining the Hilbert series, which serves as an important invariant revealing information about the poles of $h_A(t)$, number of generators of $A$,   and the Gelfand-Kirillov dimension (GK-dimension) of $A$. Indeed, the following lemma is considered well known in the study of GK dimension. 

\lemma\label{lemma3.1}(\cite{REYES2019201}, Lemma 2.7) Let $0 \neq h(t) = r(t)s(t)^{-1}$ where $r(t) \in \mathbb{Z}[t], s(t) \in \mathbb{Z}[t]$ with $s(0) = \pm1$, and $r(t), s(t)$ are relatively prime. Suppose that its Laurent series expansion $h(t) = \sum_{n \geq n_0} a_n t^n \in \mathbb{Z}((t))$ centered at zero has nonnegative coefficients. Let
\[
m = \limsup_{n \to \infty} \log_n \left( \sum_{i=n_0}^n a_i \right)
\]
denote the GK dimension of $h(t)$. Then:
\begin{enumerate}
    \item The following are equivalent:
    \begin{enumerate}
        \item $m < \infty$.
        \item The roots of $s(t)$ are all roots of unity;
        
    \end{enumerate}
    \item If $m < \infty$, then $m$ is equal to the multiplicity of $t = 1$ as a root of $s(t)$ (equivalently, the order of the pole of $h(t)$ at $t = 1$).\\
\end{enumerate}

In light of these observations,  we introduce a new class of algebras where the trivial $A$-module $\mathbb{K}_A$ has a finite free resolution characterized by monotonically increasing Betti numbers. This class not only provides partial answers to several open questions posed in the introduction but also encompasses a wide variety of algebras that have been central to ongoing research.

\definition\label{3.2} Let \( A \) be a graded algebra of finite global dimension \( d \). We say \( A \) is a \textit{monotonic} algebra if its trivial \( A \)-module \( \mathbb{K}_A \) admits a finite free resolution as in \((3.01)\), and its Betti numbers satisfy:

\[
\max \{ \alpha_{i,s} \;|\; 1 \leq s \leq n_i \} \leq \min \{ \alpha_{i+1,s} \;|\; 1 \leq s \leq n_{i+1} \}
\]
for all \( 0 \leq i \leq d \).

\medskip

We note that the existence of a finite free resolution of \( \mathbb{K}_A \) is considered a rather weak assumption since it holds for many classes of interest. For example, if \( A \) is Noetherian, every finitely generated \( A \)-module \( M \) has a finite free resolution. Similarly, Artin-Schelter Gorenstein algebras of finite global dimension satisfy this property, as shown in \cite{Stephenson1997}.

\vspace{0.5cm}

We now recall another important class of algebras, which will prove to be monotonic, those with \textit{pure} minimal resolutions.

\definition\label{3.3} Let \( A \) be a graded algebra of finite global dimension \( d \). \( A \) is said to have a \textit{pure} minimal resolution if its trivial \( A \)-module \( \mathbb{K}_A \) admits a finite free resolution as in (\ref{3.0.1}), where each \( P_i \) is generated in a single degree. In other words, \( P_i = \bigoplus_{s=1}^{n_i} A(-\alpha_i) \) for some integer \( \alpha_i \).

\vspace{0.3cm}

Several well-known classes of algebras are known to admit pure minimal resolutions. Prominent examples include Koszul algebras, \( N \)-Koszul algebras, Piecewise Koszul algebras, and \( (p, \lambda) \)-Koszul algebras. Definitions of these algebras can be found in \cite{pei2010koszul}.

\vspace{0.5cm}

We are now ready to present our first theorem of this section, which demonstrates that the class of monotonic algebras is quite extensive, encompassing numerous families of algebras that have been thoroughly studied in the literature.

\theorem\label{theorem3.4} The class of monotonic algebras includes:
\begin{enumerate}

\item[(a)] Algebras with \textit{pure} minimal resolutions.\\

\item[(b)] Artin-Schelter regular algebras of dimensions at most four.\\

\item[(c)] Quadratic Artin-Schelter regular algebras of dimensions at most six.\\

\end{enumerate}

\begin{proof} Let $A$ be a graded algebra of global dimension $d$ with a finite free resolution as in $(\ref{3.0.1})$. It follows from the minimality of  $(\ref{3.0.1})$ that 
\[  \min \;\{ \alpha_{i,s} \;\; | \;\; 1\leq s \leq n_i\}   <     \min \;\{ \alpha_{i+1,s} \;\; | \;\; 1\leq s \leq n_{i+1}\}   \tag{3.5.1}\label{3.5.1} \]
for all $0\leq i \leq d$. Furthermore, if $A$ is Gorenstein, it follows from the proof of Proposition 3.1 in  \cite{Stephenson1997} that 
\[  \max\;\{ \alpha_{i,s} \;\; | \;\; 1\leq s \leq n_i\}   <     \max \;\{ \alpha_{i+1,s} \;\; | \;\; 1\leq s \leq n_{i+1}\}   \tag{3.5.2}\label{3.5.2}\\ \]

\noindent$(a)$ Suppose that $A$ has a finite free pure resolution. Then for each $0\leq i <d$, there exists a positive integer  $\alpha_i$ such that $\alpha_i = \alpha_{i,s}$ for all  $1\leq s \leq n_i$ and so 
\[\max \{ \alpha_{i,s} \;|\; 1 \leq s \leq n_i \} = \alpha_i = \min \{ \alpha_{i,s} \;|\; 1 \leq s \leq n_i \} \]\\

Combined with inequality (\ref{3.5.1}), one has that for all $0\leq i <d$,

\[  \max \;\{ \alpha_{i,s} \;\; | \;\; 1\leq s \leq n_i\}  =  \min \;\{ \alpha_{i,s} \;\; | \;\; 1\leq s \leq n_i\} \]

\[<  \min \;\{ \alpha_{i+1,s} \;\; | \;\; 1\leq s \leq n_{i+1}\}     \]
and so $A$  is monotonic. \\

\noindent $(b)$ For the sake of brevity, we assume that $A$ is Artin-Schelter regular of dimension four since the proof is analogous in smaller dimensions. Incorporating inequalities (\ref{3.5.1}),(\ref{3.5.2}), and  Lemma \ref{lemma2.2}, it follows that the minimal free resolution of the trivial $A$-module $\mathbb{K}_A$ admits the form 

\[    A(-l) \longrightarrow A( -l+1)^{n_1} \longrightarrow  \bigoplus_{s=1}^{n_2}A(-\alpha_{2,s})   \longrightarrow  A( -1)^{n_1} \longrightarrow  A  \longrightarrow  \mathbb{K}  \]
 where $1 < \alpha_{2,s} < l-1$ for all $1\leq s \leq n_2$ and thus $A$ is monotonic.\\
 
 \noindent $(c)$ We assume that $A$ is quadratic Artin-Schelter regular of dimension six since the proof is analogous in smaller dimensions. By applying inequalities (\ref{3.5.1}), (\ref{3.5.2}), and  Lemma \ref{lemma2.2}, the minimal free resolution of the trivial $A$-module $\mathbb{K}_A$ has the form 
 
 $$  \longrightarrow \bigoplus_{i=1}^{n_2} A(-2) \longrightarrow A(-1)^{n_1}\longrightarrow A\longrightarrow  \mathbb{K}          $$
 $$    A(-l) \longrightarrow A(-l+1)^{n_1} \longrightarrow \bigoplus_{s=1}^{n_2} A(-l+2)\longrightarrow \bigoplus_{s=1}^{n_3} A(-\alpha_{i,s})    $$  
 where $2 < \alpha_{3,s} < l-2$ for all $1\leq s \leq n_3$ from which its clear that $A$ is monotonic. \end{proof}
\vspace{0.2cm}

For an example of a non-monotonic algebra, we refer the reader to Remark \ref{remark3.7} below.\\

A longstanding conjecture in the study of Artin-Schelter regular algebras, and simultaneously  in the study of Noetherian graded algebras of finite global dimension, is that the global dimension equals the GK-dimension. As a step toward this goal, the next theorem demonstrates that for monotonic algebras, the global dimension provides an upper bound for the GK-dimension.\\

 Prior to proving it, we need to recall the famous Descartes' rule of signs. 
 
 \lemma\label{lemma3.5}(Descartes' rule of signs) Given a polynomial $f(x) = \sum_{i=0}^n a_ix^{b_i}\in \mathbb{R}[x]$ with integer powers $0 \leq b_0<b_1\cdots <b_n$  and nonzero coefficients $a_i \neq0$, the number of strictly positive roots (counting multiplicity) of $f$ is equal to the number of sign changes in the coefficients of $f,$ minus a nonnegative even integer.\\

\theorem\label{theorem3.6}Let $A$ be a monotonic algebra of global dimension $d$. If $A$ has finite GK dimension, then GK$\dim(A) \leq  d$.

\begin{proof} By definition, the trivial \( A \)-module \( \mathbb{K}_A \) has a finite free resolution as in (\ref{3.0.1}). Additionally, the Hilbert series of \( A \) is given by \( h_A(t) = \frac{1}{p(t)} \), where 
\[
p(t) = \sum_{i=0}^d (-1)^i \sum_{s=0}^{n_i} t^{\alpha_{i,s}}, \tag{3.6.1}\label{3.6.1}
\]
with the \( \alpha_{i,s} \)'s and \( n_i \)'s arising from the finite free resolution. Since \( A \) is monotonic, the coefficients of \( p(t) \) change signs exactly \( d \) times. By Descartes' rule of signs, this provides an upper bound for the number of positive roots of \( p(t) \), counting multiplicities.

If \( A \) has finite GK-dimension, then \( t = 1 \) is the only positive root of \( p(t) \), as all the roots are roots of unity by Lemma \ref{lemma3.1}. Furthermore, the GK-dimension of $A$ is precisely equal to the the multiplicity of \( t = 1 \), which is bounded by $d$ through Descartes' rule of signs.
\end{proof}
\vspace{0.3cm}

\remark\label{remark3.7} Consider a graded algebra \( A \) with global dimension \( d \) whose Hilbert series takes the form \( h_A(t) = \frac{1}{p(t)} \). One might wish that the polynomial \( p(t) \) always changes signs exactly \( d \) times. However, this is not the case, indicating that some variation of the monotonicity condition may be necessary to establish such an upper bound for an even broader classes of algebras, such as Noetherian algebras of finite global dimension or all Artin-Schelter regular algebras. Indeed, in \cite{floystad2010artinschelter}, the authors provide an example of a Noetherian Artin-Schelter regular algebra of dimension five, generated by two elements, with the following Hilbert series:
\[
h_A(t) =  \frac{1}{1 -2t +t^3+t^4-t^5+t^7-t^8 -t^9+2t^{11} -t^{12}}.
\]
While \( A \) has a global dimension of five, the coefficients of \( p(t) \) exhibit seven sign changes so  $A$ is not monotonic. \\

This concludes our results regarding GK dimensions. We now turn our attention to the minimal number of generators. In \cite{lu2004regular}, the authors posed the question of whether the global dimension of an Artin-Schelter regular algebra serves as an upper bound for the minimal number of generators. We prove this for the class of monotonic algebras under two mild assumptions. To begin, we recall the following well-known lemma whose proof for instance can be found in \cite{nathanson2022}. \\

\lemma\label{lemma3.8}Let $p(t)$ be a polynomial with real coefficients and let $\nu_p$ denote the number of times its coefficients change sign. If $\alpha$ is a positive root of $p(t)$ such that $p(t) = (t-\alpha)q(t)$, then $\nu_q \leq \nu_p -1$.\\
 
 \theorem\label{theorem3.9} Let $A$ be a monotonic algebra of global dimension $d$. If GK$\dim A = d$, then $A$ is generated by  at most $d$ elements.\\
 
\begin{proof}On the contrary, suppose that $A$ is generated by $\beta$ elements where $\beta > d$. By Lemma \ref{lemma3.1}, the Hilbert series  of $A$ is of the form      $h_A(t) = \frac{1}{p(t)}$ where \[   p(t) = (1-t)^dq(t) = 1 - \beta t + \cdots     \]
Applying Lemma \ref{lemma3.8}, we deduce that \( \nu_q \leq \nu_p - d \). Since \( A \) is monotonic and \( p(t) \) has the form described in (\ref{3.6.1}), it changes signs precisely \( d \) times, implying \( \nu_p = d \). Therefore, it follows that \( \nu_q \leq d - d = 0 \).
  However after performing long division on $q(t)$, we see that 
\[  q(t) = \frac{p(t)}{(1-t)^d} = 1 -\big( \beta -d\big)t + \cdots      \]
Since $\beta -d >0$, it follows that $\nu_q \geq 1$ which is a contradiction. Therefore, we conclude that $A$ is minimally generated by at most $d$ elements.
 \end{proof}
 
 \vspace{0.3cm}
 
Lastly, we note that one can also use Theorem \ref{theorem3.6} to prove a similar result, under the hypothesis that the Hilbert series of $A$ is that of a weighted polynomial ring, which is again conjectured to be true for all Artin Schelter regular algebras, and Noetherian algebras of finite global dimension.\\

\theorem\label{theorem3.10} Let $A$ be a monotonic algebra of global dimension $d$. Furthermore, suppose that the Hilbert series of $A$ is that of weighted polynomial ring, i.e. \[ h_A(t) = \frac{1}{\prod_{i=1}^m (1-t^i)^{n_i}}\]
 Then $A$ is generated by at most $d$ elements.\\
 
\begin{proof} By the Hilbert series assumption, 
 \[ p(t) = \prod_{i=1}^m (1-t^i)^{n_i}= 1 -{n_1} t + \cdots \]
 By Theorem \ref{theorem3.6}, GK$\dim(A) \leq d$ and thus $n_1 \leq d$. The result is now immediate from the well known fact that $n_1$ represents the number of generators of $A$.
 
 \end{proof}

 \section*{Acknowledgments}
The author would like to extend sincere gratitude to his advisor, Manuel L. Reyes, for his indispensable mentorship, insightful feedback, and careful proofreading of this paper. This paper constitutes a portion of the author's Ph.D. thesis.\\

\printbibliography

\end{document}